\documentclass{svjour3}
\usepackage{amsfonts}
\usepackage{amsmath}

\smartqed \journalname{Acta Applicandae Mathematicae}
\begin{document}

\title{Feedback Differential Invariants}

\author{Valentin Lychagin}

\institute{V.Lychagin \at
              Department of Mathematics and Statistics, University of Tromso \\90377, Tromso, Norway\\
              \email{lychagin@math.uit.no}  \\
}

\date{Received:  / Accepted: }
\maketitle


\begin{abstract}
 The problem of feedback equivalence for control systems is considered.
 An algebra of differential invariants and criteria for the feedback equivalence for regular control systems are found.

\end{abstract}


\section{Introduction}
In this paper we outline an application of the  method of differential invariants to the problem of recognition for control systems.
Namely, we consider the action of feedback transformations on $1$-dimensional control autonomous systems describing by the second order ordinary differential equations.
It is easy to check that for these systems (and it is the general case) one has infinite number functionally independent invariants,  but, and it is also common for such problems, they could be organized in algebra of functions on some differential equation (so-called syzygy). In other words, there is a finite number of differential invariants and invariant derivations such that any differential invariant can be obtained by computing functions of these basic invariants and their derivations. This is the Lie-Tresse theorem (see, for example,\cite{T}\cite{Kum},\cite{KL},\cite{OL} ).
 In our case, for the description of the algebra in neighborhoods of regular orbits (Theorem?), we need one differential invariant of the 1st order and two differential invariants of the 3rd order. By the definition, differential invariants describe the orbits of jet of the system. Hence, the structure of the differential invariant algebra allows us to establish formal feedback equivalence of the control systems. In section ? we show that for the feedback pseudogroup the formal equivalence implies the  smooth in the case of regular systems.
 There are several approaches to study control systems. The  more popular based either on EDS methods (\cite{C},\cite{GS},\cite{H}), or on  affine families of vector fields (\cite{AZ},\cite{J},\cite{WR}). 
 In this paper we consider control systems as underdetermined differential equations. The corresponding geometrical picture leads us to submanifolds in the jet spaces and to the pseudogroup of Lie transformations. It makes very natural to consider the feedback transformations as well as produces differential invariants for problems where the EDS method did not find them. Say problems investigated in this paper are equivalent (in general case) as EDS systems, but not as control ones.

\section{Feedback Pseudogroup}

Let
\begin{equation}
y^{\prime\prime}=F(y,y^{\prime},u), \label{Cont1}%
\end{equation}
be an autonomous 1-dimensional control system.

Here the function $y=y\left(  t\right)  $ describes a dynamic of the state of the system,
and $u=u\left(  t\right)$ is a scalar control parameter.

We shall consider this system as an undetermined ordinary differential equation of the second order on sections of 
 $2$-dimensional bundle $\pi:\mathbb{R}^{3}\rightarrow\mathbb{R}$ , where
$\pi:(u,y,t)\longmapsto t.$

Let $\mathcal{E}\subset J^{2}\left(  \pi\right)  $ be the corresponding
submanifold. In the canonical jet coordinates $\left(  t,u,y,u_{1}%
,y_{1},....\right)  $ this submanifold is given by the equation:%
\[
y_{2}=F\left(  y,y_{1},u\right)  .
\]

It is known (see, for example, \cite{KLV}) that Lie transformations in jet bundles
$J^{k}\left(  \pi\right)  $ for $2$-dimensional bundle $\pi$ are prolongations
of point transformations, that is, prolongations of diffeomorphisms of the
total space of the bundle $\pi.$

We shall restrict ourselves by point transformations which are automorphisms of the
bundle $\pi.$

Moreover, if these transformations preserve the class of systems (\ref{Cont1})
then they should have the form
\begin{equation}
\Phi:\left(  u,y,t\right)  \rightarrow\left(  U\left(  u,y\right)  ,Y\left(
y\right)  ,t\right)  . \label{Feed1}%
\end{equation}

Diffeomorphisms of form (\ref{Feed1}) is called \emph{feedback
transformations. }
The corresponding infinitesimal version of this notion is a
\emph{feedback vector field, }i.e. a plane vector field of the form%
\[
X_{a,b}=a\left(  y\right)  \partial_{y}+b\left(  u,y\right)  \partial_{u}.
\]

The feedback transformations in a natural way act on the control systems of type (\ref{Cont1}):
\[
\mathcal{E}\longmapsto\Phi^{\left(  2\right)  }\left(  \mathcal{E}\right)  ,
\]
where $\Phi^{\left(  2\right)  }:J^{2}\left(  \pi\right)  \rightarrow
J^{2}\left(  \pi\right)  $ is the second prolongation of the point
transformation $\Phi.$

Passing to functions $F,$ defining the systems, we get the following
action on these functions:%
\begin{equation}
\widehat{\Phi}:F\left(  y,y_{1},u\right)  \longmapsto\frac{1}{Y^{\prime}%
}F\left(  Y,Y^{\prime}y_{1},U\right)  -\frac{Y^{\prime\prime}}{Y^{\prime}%
}y_{1}^{2}. \label{rep1}%
\end{equation}

The infinitesimal version of this action leads us to the following presentation of
feedback vector fields:%

\begin{equation}
\widehat{X_{a,b}}~=b\left(  u,y\right)  \partial_{u}+a\left(  y\right)
\partial_{y}+a^{\prime}y_{1}\partial_{y_{1}}+\left(  a^{\prime\prime}y_{1}%
^{2}+a^{\prime}f\right)  \partial_{f}. \label{rep2}%
\end{equation}

In this formula $\widehat{X_{a,b}}$ is a vector field on the $4$-dimensional
space $\mathbb{R}^{4}$ with coordinates $\left(  y,y_{1},u,f\right)  .$
Each control system (\ref{Cont1}) determines a $3$-dimensional submanifold $L_{F}%
\subset\mathbb{R}^{4},$ the graph of $F:$%
\[
L_{F}=\left\{  f=F\left(  y,y_{1},u\right)  \right\}  .
\]
Let $A_{t}$ be the $1$-parameter group
of shifts along vector field $X_{a,b}$ and let $B_{t}:\mathbb{R}^{4}\rightarrow
\mathbb{R}^{4}$ be the corresponding $1$-parameter group of shifts along
$\widehat{X_{a,b}},$ then these two actions related as follows
\[
L_{\widehat{A_{t}}\left(  F\right)  }=B_{t}\left(  L_{F}\right)  .
\]
In other words, if we consider an $1$-dimensional bundle
$$
\kappa : \mathbb{R}^4\rightarrow\mathbb{R}^3,
$$
where
$
\kappa ((  y,y_{1},u,f))= ( y,y_{1},u),
$
then formula (\ref{rep2}) defines the representation $X\longmapsto\widehat{X}$ of the
Lie algebra of feedback vector fields into the Lie algebra of Lie vector
fields on $J^{0}\left(  \kappa\right)  ,$ and the action of Lie vector
fields $\widehat{X}$ on sections of bundle $\kappa$ corresponds to the
action of feedback vector fields on right hand sides of (\ref{Cont1})(see,\cite{KLV}).

\section{Feedback Differential Invariants}

By a \emph{feedback differential invariant} of order $\leq k$ we understand a function
$I\in C^{\infty}\left(  J^{k}\kappa\right)  $ on the space of $k$-jets
$J^{k}(\kappa),$ which is invariant under of the prolonged action of
feedback transformations.

Namely,
\[
\widehat{X_{a,b}}^{\left(  k\right)  }\left(  I\right)  =0,
\]
for all feedback vector fields $X_{a,b}.$

In what follows we shall omit subscript of order of jet spaces,
and say that a function $I$ on the space of infinite jets $I\in C^{\infty
}\left(  J^{\infty}\kappa\right)  $ is a feedback differential invariant
if
$$
\widehat{X_{a,b}}^{\left(  \infty\right)  }\left(  I\right)  =0.
$$

In a similar way one defines a \emph{feedback invariant derivations} as combinations of
total derivatives%
\begin{align*}
&  \nabla=A\frac{d}{dy}+B\frac{d}{dy_{1}}+C\frac{d}{du},\\
&  A,B,C\in C^{\infty}\left(  J^{\infty}\kappa\right)  ,
\end{align*}
which are invariant with respect to prolongations of feedback transformations,
that is,%
\[
\lbrack\widehat{X_{a,b}}^{\left(  \infty\right)  },\nabla]=0
\]
for all feedback vector fields $X_{a,b}.$

Remark that for such derivations functions $\nabla\left(  I\right)  $ are differential
invariants (of order higher then order of $I$) for any feedback differential
invariant $I.$ This observation allows us to construct new differential
invariants from known ones by the differentiations only.

\section{Dimensions of Orbits}

First of all, we remark that the submanifold $y_{1}^{(-1)}(0)$ is a singular orbit for the feedback action in the space of $0$-jets  $J^{0}\kappa $. In what follows we shall consider orbits of jets at regular points, that is, at such points, where $y_{1}\neq 0$.

It is easy to see, that the $k$th prolongation of the feedback vector field $\widehat{X_{a,b}}$
depends on $\left(  k+2\right)  $-jet of function $a\left(  y\right)  $ and
$k$-jet of function $b\left(  u,y\right)  .$

Denote by $V_{i}^{k}$ and $W_{ij}^{k}$ the components of the decomposition%
\[
\widehat{X_{a,b}}^{\left(  k\right)  }=%
{\displaystyle\sum\limits_{i=0}^{k+2}}
a^{\left(  i\right)  }\left(  y\right)  V_{i}^{k}+%
{\displaystyle\sum\limits_{0\leq i+j\leq k}}
\frac{\partial^{i+j}b}{\partial u^{i}\partial y^{j}}W_{ij}^{k}.
\]
Then, by the construction, the vector fields $V_{i}^{k},0 \leq i \leq k+2$, and $W_{ij}^{k},0\leq i+j \leq k$, generate 
a completely integrable distribution on the space of $k$-jets, integral manifolds of which are orbits of the feedback action in $J^{k}\kappa.$ 

Let
$\mathcal{O}_{k+1}$ be an orbit in $J^{k+1}\kappa$, then the projection
$\mathcal{O}_{k}=\kappa_{k+1,k}\left(  \mathcal{O}_{k+1}\right)  \subset
J^{k}\kappa$ is an orbit too, and to determine dimensions of the orbits one
should find dimensions of the bundles: $\kappa_{k+1,k}:\mathcal{O}%
_{k+1}\rightarrow\mathcal{O}_{k}.$ To do this we should find conditions on
functions $a$ and $b$ under which $\widehat{X_{a,b}}^{\left(  k\right)  }=0$
at a point $x_{k}\in J^{k}\kappa.$

It easy to check that these conditions for $k=1$ at a point $x_{1},$ where
$f_{u}\neq0,$ has the form:%
\begin{align}
a\left(  y\right)   &  =\ a^{\prime}\left(  y\right)  =\ a^{\prime\prime
}\left(  y\right)  =0,b\left(  u,y\right)  =0,\label{sing1}\\
b_{u}  &  =0,a^{\prime\prime\prime}y_{1}^{2}-b_{y}f_{u}=0.\nonumber
\end{align}
Here $\left(  u,y,y_{1},f,f_{u},f_{y},f_{y_{1}}\right)  $ are the canonical
coordinates in the $1$-jet space $J^{1}(\kappa).$

The formula for prolongations of vector fields (see, for example,\cite{KLV}) shows that 
the conditions on functions $a$ and $b$ such that vector fields $\widehat{X_{a,b}}^{\left(
k\right)  }$ vanish at a point in $J^{k}\kappa$ are just $\left(  k-1\right)  $-prolongations
of (\ref{sing1}).

Let%
\[
\phi=a^{\prime\prime}\left(  y\right)  y_{1}^{2}+a^{\prime}\left(  y\right)
f-b\left(  u,y\right)  f_{u}-a\left(  y\right)  f_{y}-a^{\prime}\left(
y\right)  y_{1}f_{y_{1}}%
\]
be the generating function of vector field $\widehat{X_{a,b}}.$

Assume that $k>1,$ and that $\widehat{X_{a,b}}^{\left(  k-1\right)  }=0$ at a
point $x_{k-1}\in J^{k-1}\kappa$ . Then the vector field $\widehat{X_{a,b}%
}^{\left(  k\right)  }$ is a $\kappa_{k,k-1}$-vertical over this point.
Components
\[
\frac{d^{k}\phi}{du^{i}dy^{j}}\frac{\partial}{\partial f_{\sigma_{ij}}}%
\]
of this vector field, where $\sigma_{ij}=(\underset{i\text{-times}%
}{\underbrace{u,....,u},}\underset{j\text{-times}}{\underbrace{y...,y}}), i+j=k,$
and components
\[
\frac{d^{k}\phi}{dy^{k-1}dy_{1}}\frac{\partial}{\partial f_{\tau}},
\]
where $\tau=(\underset{(k-1)\text{-times}}{\underbrace{y...,y}},y_{1})$
depend on
\[
\frac{\partial^{k}b}{\partial u^{i}\partial y^{j}}%
\]
and
\[
\frac{d^{k+2}a}{dy^{k+2}}%
\]
respectively.

All others components
\[
\frac{d^{k}\phi}{dy^{r}dy_{1}^{s}du^{t}}\frac{\partial}{\partial f_{\sigma}}%
\]
are expressed in terms of $\left(  k-1\right)  $-jet of $b\left(  u,y\right)
$ and $\left(  k+1\right)  $-jet of function $a\left(  y\right)  .$

It shows that the bundles: $\kappa_{k,k-1}:\mathcal{O}_{k}\rightarrow
\mathcal{O}_{k-1}$ are $(k+2)$-dimensional if $k>1$, and $y_{1}\neq0,f_{u}%
\neq0.$

 We say that $k$-jet $[F]^{k}_{p}\in J^{k}\kappa$ of a function $F$ is \emph{weakly regular} if the point $p$ is regular, that is, $y_{1}\neq 0$ at this point, and $F_{u}\neq 0$.

  Orbits of the weakly regular points we call \emph{weakly regular}.

  Feedback orbits in the space of $1$-jets can be found by direct integration of $6$-dimensional completely integrable distribution generating by the vector fields $V^{1}_{i}, 0\leq i \leq 3$, and $W^{1}_{ij}, 0\leq i+j \leq 1$.
  Summarizing, we get the following result.
\begin{theorem}
\begin{enumerate}
\item The first non-trivial differential invariants of feedback
transformations appear in order $1$ and they are functions of the basic
invariant
\[
J=\frac{y_{1}f_{y_{1}}-2f}{y_{1}}.
\]

\item Dimension of weakly regular orbit of feedback transformations in $J^{k}%
\kappa,$ $k>1,$ is equal to%
\[
\frac{\left(  k+2\right)  \left(  k+3\right)  }{2}.
\]

\item There are
\[
\frac{\left(  k+2\right)  \left(  k-1\right)  }{2}%
\]
independent differential invariants of pure order $k.$
\end{enumerate}
\end{theorem}

\section{Invariant Derivations}

We expect three linear independent feedback invariant derivations. The straightforward computations in order $\leq2$ show that they are of the form
\begin{align*}
\nabla_{u}  &  =\frac{y_{1}}{f_{u}}\frac{d}{du},\\
\nabla_{y}  &  =-{\frac{{y}_{1}^{3}f_{{y}_{1}u}-2\,{z}^{2}f_{{y}}+{y}_{1}%
^{2}f~f_{{y}_{1}y_{1}}-2y_{1}\,f\ f_{{y}_{1}}+2\,{f}^{2}}{y_{1}\left(
-2\,f_{{u}}+y_{1}f_{{y}_{1}{u}}\right)  }}\frac{d}{du}+y_{1}\frac{d}%
{dy}+f\frac{d}{dy_{1}},\\
\nabla_{y_{1}}  &  =y_{1}\frac{d}{dy_{1}}.
\end{align*}

It is easy to check that these derivations obey the following commutation relations

\begin{eqnarray}
\left[ \nabla _{u},\nabla _{y}\right]  &=&\frac{L-J_{y_{1}y_{1}}}{J_{u}}%
\nabla _{u}+\nabla _{y_{1}},  \label{CommutationRelations} \\
\left[ \nabla _{u},\nabla _{y_{1}}\right]  &=&\left( 1+J_{u}\right) ~\nabla
_{u},  \notag \\
\left[ \nabla _{y},\nabla _{y_{1}}\right]  &=&\frac{J_{u}K+J_{y_{1}}\left(
J_{y_{1}}-J_{u}+J~J_{u}\right) -J_{y_{1}y_{1}}\left( J_{y_{1}}-J_{u}\right)
}{J_{u}^{2}}~\nabla _{u}  \notag \\
&&-\nabla _{y}-J~\nabla _{y_{1}},  \notag
\end{eqnarray}
where $K$ and $L$ are differential invariants of the 3rd order (see below).

\section{Differential Invariants of Order $2$}

Theorem 1 shows that there are $2$ independent differential invariants of pure order $2.$
We can get them by applying invariant derivations to the $1$st order invariant $J$:%
\begin{align*}
J_{u}\overset{\text{def}}{=}\nabla_{u}\left(  J\right)   &  =\frac
{y_{1}f_{y_{1}u}-2f_{u}}{f},\\
J_{y_{1}}\overset{\text{def}}{=}\nabla_{y_{1}}\left(  J\right)   &
=\frac{y_{1}^{2}f_{y_{1}y_{1}}-2y_{1}f_{y_{1}}+2f}{y_{1}^{2}},
\end{align*}

but

\[
\nabla_{y}\left(  J\right)  =0.
\]

\section{Differential Invariants of Order $3$}

Theorem 1 shows that there are five independent differential invariants of the $3$rd order. Three of them
we get by invariant differentiation:%
\[
J_{u~u}=\nabla_{u}\nabla_{u}\left(  J\right)  ,\ J_{u~y_{1}}=\nabla_{u}%
\nabla_{y_{1}}\left(  J\right)  ,\ J_{y_{1}y_{1}}=\nabla_{y_{1}}\nabla_{y_{1}%
}\left(  J\right)  .
\]

To find the last 2 differential invariants  we remark that the $3$-prolongations of feedback vector fields are affine along the fibres
$$
J^{3}\left( \kappa \right) \overset{\kappa _{3,2}}{\rightarrow }J^{2}\left(
\kappa \right)
$$
see,for example, (\cite{KLV}).

Therefore one can try to find the differential invariants as functions which are affine along the fibres $\kappa_{3,2}$.

Finally, we get:

\begin{align*}
K  &  =y_{1}^{2}f_{uy_{1}y_{1}}-\frac{3y_{1}~J_{u}~-y_{1}J_{y_{1}u}}{J_{u}%
}f_{yy_{1}}-y_{1}J_{u}f_{yy_{1}}+\frac{2(J_{y_{1}u}+2J_{u})}{J_{u}}f_{u}\\
&  +2J_{u}f_{y}-\frac{(J_{u}J_{y_{1}}+J_{y_{1}}-J_{y_{1}y_{1}})J_{u}+J_{y_{1}%
}J_{y_{1}u}}{y_{1}\ J_{u}}f+\frac{(J_{u}-J_{y_{1}})(J_{y_{1}y_{1}}+J_{y_{1}}%
)}{J_{u}},
\end{align*}

and%

\begin{align*}
L  &  =\frac{y_{1}^{2}}{f_{u}}f_{uyy_{1}}-\frac{(2+J_{u})~y_{1}}{f_{u}}%
f_{uy}-\frac{y_{1}J_{uu}}{J_{u}}f_{yy_{1}}+2\frac{J_{uu}}{J_{u}}f_{y}+\left(
\frac{J_{y_{1}u}}{y_{1}}-\frac{J_{y_{1}}J_{uu}}{y_{1}J_{u}}\right)  f\\
&  +J_{y_{1}}+J_{y_{1}y_{1}}.
\end{align*}

\section{Algebra of Feedback Differential Invariants}
To use the above computations one should reinforce the notion of regularity. We give the following definition.
\begin{definition}
We say that a weakly regular orbit is regular if $J_{u}\neq0,$ on the orbit.
\end{definition}

Remark that for \emph{irregular} or \emph{singular} control systems one has $J_{u}\equiv0$, and therefore they have the form:
$$
y^{\prime\prime}=A(y,y^{\prime})+B(u,y){y^{\prime}}^{2}.
$$

Counting dimensions shows that differential invariants $J,K,L$ are generators in the algebra of feedback differential invariants, and considering symbols of differential invariants shows that they satisfy two syzygy relations.

\begin{theorem}
\label{MainDiffInv}

\begin{enumerate}
\item Algebra of feedback differential invariants in a neighborhood of
regular orbits is generated by differential invariant $J$ of the $1$-st order,
differential invariants $K$ and $L$ of the $3$-rd order and all their
invariant derivatives.

\item Syzygies for this algebra have two generators of the form%
\begin{align*}
J_{y} &  =0,\\
K_{u}-L_{y_{1}}+\frac{J_{y_{1}u}+J_{u}-J_{x}^{2}}{J_{u}}L-\frac{J_{u~u}}%
{J_{u}}K &  =\Phi\left(  J,J_{u},J_{y_{1}}\right)  .
\end{align*}

\end{enumerate}
\end{theorem}

\bigskip

\begin{remark}
In a similar way, for irregular systems, we get the following description of
differential invariants algebra.

Algebra\ of differential invariants for systems with $J_{u}\equiv0,$ but
$y_{1}\neq0,$ is generated by differential invariant $J$ of the $1$-st order,
differential invariant $M$ of the $3$-rd order%
\[
M=y_{1}f_{y_{1}y_{1}y_{1}}+y_{1}^{2}f_{y~y_{1}y_{1}}-ff_{y_{1}y_{1}~}%
-2y_{1}f_{y_{1}y~}+\frac{2f~f_{y_{1}}}{y_{1}}+2f_{y}-\frac{2f^{2}}{y_{1}^{2}}%
\]
and all invariant derivatives
\[
\nabla_{y_{1}}^{a}J,
\]%
\[
\nabla_{y_{1}}^{a}\nabla_{u}^{b}M.
\]

\end{remark}

\section{The Feedback Equivalence Problem}

Consider two control systems given by functions $F$ and $G.$ Then, to establish
feedback equivalence, one should solve the differential equation%
\begin{equation}
\frac{1}{Y^{\prime}}F\left(  Y,Y^{\prime}y_{1},U\right)  -\frac{Y^{\prime
\prime}}{Y^{\prime}}y_{1}^{2}-G\left(  y,y_{1},u\right)  =0\label{equiveq1}%
\end{equation}
with respect to functions $Y\left(  y\right)  $ and $U\left(  u,y\right)  .$

Let us denote the left hand side of (\ref{equiveq1}) by $H.$ Then assuming the
general position one can find functions $U,Y,Y^{\prime},Y^{\prime\prime}$ from
the equations%
\[
H=H_{y_{1}}=H_{y_{1}y_{1}}=H_{y_{1}y_{1}y_{1}}=0.
\]
Assume that we get%
\begin{align*}
U  & =A\left(  y,y_{1},u\right)  ,Y=B\left(  y,y_{1},u\right)  ,\\
Y^{\prime}  & =C\left(  y,y_{1},u\right)  ,Y^{\prime\prime}=D\left(
y,y_{1},u\right)  .
\end{align*}
Then the conditions
\begin{align*}
A_{y_{1}}  & =B_{y_{1}}=C_{y_{1}}=D_{y_{1}}=0,\\
B_{u}  & =C_{u}=D_{u}=0
\end{align*}
and
\[
C=B_{y},D=C_{y}%
\]
show that if (\ref{equiveq1}) has a formal solution at each point $\left(
y,y_{1},u\right)  $ in some domain then this equation has a smooth solution.

On the other hand if system $F$ at a point $p=(y^{0},y_{1}^{0},u^{0})$ and
system $G$ at a point $\widetilde{p}=(\widetilde{y}^{0},\widetilde{y}_{1}%
^{0},\widetilde{u}^{0})$ has the same differential invariants then, by the
definition, there is a formal feedback transformation which send the infinite
jet of $F$ at the point $p$ to the infinite jet of $G$ at the point
$\widetilde{p}.$

 Keeping in mind these observations and results of theorem
\ref{MainDiffInv} we consider the space $\mathbb{R}^{3}$with coordinates $\left( u,y,y_{1}\right)
$ and the space $\mathbb{R}^{8}$ with coordinates $\left(
j,j_{1},j_{3},j_{11},j_{13},j_{33},k,l\right) .$

Then any control system,
given by the function $F\left( u,y,y_{1}\right) $, defines a map
\begin{equation*}
\sigma _{F}:\mathbb{R}^{3}\rightarrow \mathbb{R}^{8},
\end{equation*}%
by
\begin{eqnarray*}
j &=&J^{F},j_{1}=J_{u}^{F},j_{3}=J_{y_{1}}^{F}, \\
j_{11} &=&J_{u~u}^{F},j_{13}=J_{u~y_{1}}^{F},j_{33}=J_{y_{1}y_{1}}^{F}, \\
k &=&K^{F},l=L^{F},
\end{eqnarray*}
where the subscript $F$ means that the differential invariants are evaluated
due to the system.

Let
\begin{equation*}
\Phi :\mathbb{R}^{3}\rightarrow \mathbb{R}^{3}
\end{equation*}%
be a feedback transformation.

Then from the definition of the feedback
differential invariants it follows that
\begin{equation*}
\sigma _{F}\circ \Phi =\sigma _{\widehat{\Phi }\left( F\right) }.
\end{equation*}%
Therefore, the geometrical image
\begin{equation*}
\Sigma _{F}={Im}\left( \sigma _{F}\right)\subset\mathbb{R}^{8}
\end{equation*}%
does depend on the feedback equivalence class of $F$ only.

We say that a system $F$ is \emph{regular} in a domain $D\subset\mathbb{R}%
^{3}$ if

\begin{enumerate}
\item $3$-jets of $F$ belong to regular orbits,

\item  $\sigma_{F}\left(  D\right)  $ is a smooth $3$-dimensional submanifold in $\mathbb{R}%
^{8},$ and

\item functions $j,j_{1},j_{3}$ are coordinates on $\Sigma
_{F}.$
\end{enumerate}

The following lemma gives a relation between the Tresse derivatives and invariant differentiations $\nabla _{u},\nabla _{y},\nabla _{y_{1}}$.
\begin{lemma}
Let
$$
\frac{D}{DJ},\frac{D}{DJ_{u}},\frac{D}{DJ_{y_{1}}}
$$
be the Tresse
derivatives with respect to differential invariants $J,J_{u}$ and $J_{y_{1}}.$

Then the following decomposition
\begin{eqnarray*}
\nabla _{u} &=&J_{u}\frac{D}{DJ}+J_{u~u}\frac{D}{DJ_{u}}+J_{u~y_{1}}\frac{D}{%
DJ_{y_{1}}}, \\
\nabla _{y} &=&\left( J_{y_{1}y_{1}}-J_{y_{1}}-L\right) \frac{D}{DJ_{u}}+%
\frac{J_{u}K+J_{y_{1}}\left( J_{y_{1}}-J_{u}\right) -J_{y_{1}y_{1}}\left(
J_{y_{1}}-J_{u}\right) }{J_{u}}\frac{D}{DJ_{y_{1}}}, \\
\nabla _{y_{1}} &=&J_{y_{1}}\frac{D}{DJ}+\left(
J_{u~y_{1}}-J_{u}-J_{u}^{2}\right) \frac{D}{DJ_{u}}+J_{y_{1}y_{1}}\frac{D}{%
DJ_{y_{1}}}
\end{eqnarray*}%
holds.
\end{lemma}

\begin{proof}
The proof follows directly from the definition of the Tresse derivatives and
commutation relations (\ref{CommutationRelations}).
\end{proof}
\begin{theorem}
Two regular systems $F$ and $G$ are locally feedback equivalent if and only
if
\begin{equation}
\Sigma_{F}=\Sigma_{G}.\label{EquivCond}%
\end{equation}

\end{theorem}

\begin{proof}
Let us show that the condition \ref{EquivCond} implies a local feedback
equivalence.

Assume that
\begin{align*}
J_{uu}^{F}  & =j_{11}^{F}\left(  J^{F},J_{u}^{F},J_{y_{1}}^{F}\right)
,J_{uy_{1}}^{F}=j_{12}^{F}\left(  J^{F},J_{u}^{F},J_{y_{1}}^{F}\right)
,J_{y_{1}y_{1}}^{F}=j_{22}^{F}\left(  J^{F},J_{u}^{F},J_{y_{1}}^{F}\right)
,\\
K^{F}  & =k^{F}\left(  J^{F},J_{u}^{F},J_{y_{1}}^{F}\right)  ,L^{F}%
=l^{F}\left(  J^{F},J_{u}^{F},J_{y_{1}}^{F}\right)
\end{align*}
on $\Sigma_{F},$ and
\begin{align*}
J_{uu}^{G}  & =j_{11}^{G}\left(  J^{G},J_{u}^{G},J_{y_{1}}^{G}\right)
,J_{uy_{1}}^{G}=j_{12}^{G}\left(  J^{G},J_{u}^{G},J_{y_{1}}^{G}\right)
,J_{y_{1}y_{1}}^{G}=j_{22}^{G}\left(  J^{G},J_{u}^{G},J_{y_{1}}^{G}\right)
,\\
K^{G}  & =k^{G}\left(  J^{G},J_{u}^{G},J_{y_{1}}^{G}\right)  ,L^{G}%
=l^{G}\left(  J^{G},J_{u}^{G},J_{y_{1}}^{G}\right)
\end{align*}
on $\Sigma_{G}.$

Then condition  \ref{EquivCond} shows that  $j_{11}^{F}=j_{11}^{G},j_{12}%
^{F}=j_{12}^{G},j_{22}^{F}=j_{22}^{G}$ and $k^{F}=k^{G},l^{F}=l^{G}.$
Moreover,as we have seen the invariant derivations $\nabla_{u},\nabla_{y},\nabla_{y_{1}}$ are
linear combinations of the Tresse derivatives.

In other words, functions
$$
j_{11}^{F},j_{12}^{F},j_{22}^{F},k^{F},l^{F}
$$
and their partial derivatives in $j,j_{1},j_{3}$ determine the restrictions of all differential invariants.

Therefore, condition \ref{EquivCond} equalize  restrictions of differential invariants
not only to order $\leq3$ but in all orders, and provides therefore formal feedback
transformation between $F$ and $G$.

 The resulting feedback transformation has
the form%
\[
D_{F}\overset{\kappa_{F}}{\rightarrow}{\Sigma}_{F}=\Sigma_{G}%
\overset{\kappa_{G}^{-1}}{\rightarrow}D_{G},
\]
where $D_{F}$ and $D_{G}$ are domains of definition for system $F$ and $G$
respectively. 

Remark that $\kappa_{G}^{-1}\circ\kappa_{F}$ is a feedback
transformation because it sends trajectories of vector field $\nabla_{u}$ to
themselves, that are fibres of the bundle $\pi.$
\end{proof}


\end{document}